\documentclass{amsart}

\usepackage{times,amssymb,booktabs,amsmath,tikz,wasysym,stmaryrd}
\usepackage[mathscr]{euscript}
\usepackage{cite}

\usepackage[a4paper=true,pagebackref=true]{hyperref}
\renewcommand*{\backref}[1]{}
\renewcommand*{\backrefalt}[4]{{\tiny[%
  \ifcase #1 Not cited.\relax\or Page~#2%
  \else Pages #2\fi]}}
  \hypersetup{pdftitle={Blocks of the truncated $q$-Schur algebras of type A},
  pdfauthor={Andrew Mathas and Marcos Soriano},
  colorlinks = true,
  linkcolor =blue,
  anchorcolor = red,
  citecolor = blue,
  pagecolor = red,
  urlcolor = blue
}

\let\Sect=\S

\numberwithin{equation}{section}

\usepackage{aliascnt}
\def\NewTheorem#1{%
  \newaliascnt{#1}{equation}
  \newtheorem{#1}[#1]{#1}
  \aliascntresetthe{#1}
  \expandafter\def\csname #1autorefname\endcsname{#1}
}
\def\equationautorefname~#1\null{(#1)\null}

\newcounter{main}
\theoremstyle{plain}

\newtheorem*{MainTheorem}{Main Theorem}

\NewTheorem{Assumption}
\NewTheorem{Proposition}
\NewTheorem{Theorem}
\NewTheorem{Corollary}

\NewTheorem{Lemma}
\NewTheorem{Definition}
\theoremstyle{definition}
\theoremstyle{remark}
\NewTheorem{Remark}

\newaliascnt{Example}{equation}
\aliascntresetthe{Example}

\newenvironment{Example}%
  {\refstepcounter{Example}\trivlist
   \item[\hskip\labelsep\theExample.~\textbf{Example}\space]
   \ignorespaces
  }{\unskip\nobreak\hfil%
    \penalty50\hskip2em\hbox{}\nobreak\hfil$\Diamond$%
    \parfillskip=0pt\finalhyphendemerits=0\penalty-100\endtrivlist
    \medskip
}

\DeclareMathAlphabet{\mathpzc}{OT1}{pzc}{m}{it}

\let\gdom=\triangleright

\let\gedom=\trianglerighteq
\let\<\langle
\let\>\rangle

\def\Diag(#1){\llbracket#1\rrbracket}

\def\pmod#1{\hspace*{0.7ex}(\mathop{\rm mod }#1)}

\def\({\big(}
\def\){\big)}

\def\Mod{\text{Mod-}}

\def\Sym{\mathfrak S}
\def\H{\mathscr{H}}
\def\N{\mathbb N}
\def\p{\mathfrak p}
\def\P{\mathpzc{P}_{e,p}}
\def\O{\mathscr O}
\def\S{\mathscr{S}}
\def\C{\mathbb C}
\def\Z{\mathbb Z}

\def\map#1#2{\,{:}\,#1\!\longrightarrow\!#2}

{\catcode`\|=\active
  \gdef\set#1{\mathinner{\lbrace\,{\mathcode`\|"8000%
                                   \let|\midvert #1}\,\rbrace}}
}
\def\midvert{\egroup\mid\bgroup}

\newcommand\core[1][e]{\mathop{\text{core}}_{#1}\nolimits}
\DeclareMathOperator{\End}{End}

\DeclareMathOperator{\res}{res}
\DeclareMathOperator{\rad}{rad}

\makeatletter
\def\abacus#1#2#3{
  \begin{tikzpicture}[scale=0.4]
      \@tempcnta=#1\advance\@tempcnta by -1
      \@tempcntb=#2\advance\@tempcntb by -1
      \foreach \x in {0,1,...,\the\@tempcnta} {
        \draw[thick] (\x,0.3)--(\x,-\the\@tempcntb)--++(0,-0.3);
	\foreach \y in {0,1,...,\the\@tempcntb} {
	  \draw[thin,shift={(\x,-\y)}](-0.2,0)--(0.2,0);
	}
      }
      \foreach \z in { #3 } {
        \@tempcnta=\z\divide\@tempcnta by#1
        \@tempcntb=\@tempcnta\multiply\@tempcntb by #1
        \advance\@tempcntb by -\z
        \multiply\@tempcnta by -1
        \multiply\@tempcntb by -1
        \fill(\the\@tempcntb,\the\@tempcnta) circle (6pt);
      }
  \end{tikzpicture}
}
\makeatother

\def\Email#1{\email{\href{mailto:#1}{#1}}}
\begin{document}
\bibliographystyle{andrew}
\title{Blocks of the truncated $q$-Schur algebras of type A}
\author{Andrew Mathas}
\address{School of Mathematics and Statistics F07, University of
Sydney, NSW 2006, Australia.}
\Email{andrew.mathas@sydney.edu.au}
\author{Marcos Soriano}
\address{Institut f\"ur Algebra, Zahlentheorie und Diskrete Mathematik,
         Im Welfengarten 1, Leibniz Universit\"at Hannover, Deutschland}
\Email{soriano@math.uni-hannover.de}

\begin{abstract} This paper classifies the blocks of the truncated $q$-Schur algebras
  of type~$A$ which have as weight poset an arbitrary cosaturated set of partitions.
\end{abstract}

\maketitle

\section{Introduction}

In this paper we classify the blocks of the truncated $q$-Schur algebras of type
$A$. The truncated Schur algebras are a natural family of quasi-hereditary
algebras obtained from the $q$-Schur algebras by applying ``Schur functors''.
These algebras include, as special cases, the usual Schur algebras
$\S_{k,q}(n,r)$. Our main result gives a relatively quick and easy
classification of the blocks of all of these algebras.

We think it quite remarkable that there is a uniform and relatively simple
classification of the blocks of all of the truncated $q$-Schur algebras. As with
the original classification of the blocks of the $q$-Schur algebra
$\S_{k,q}(n,n)$~\cite{JM:schaper}, the main tool that we use is the Jantzen sum
formula, however, the new theme which permeates this paper is that the
\textit{Jantzen coefficients}, the integers which appear in the Jantzen sum
formula, determine much of the representation theory. For example, the blocks
correspond to the combinatorial ``linkage classes'' of partitions
determined by the non-zero Jantzen coefficients.

The proof of the classification of blocks of the Schur algebras that we give is
new, both for $\S_{k,q}(n,n)$ and more generally for the algebras
$\S_{k,q}(n,r)$ considered in \cite{Donkin:SchurIV,Cox:blocks}. Throughout our
focus is on the combinatorics of the Jantzen coefficients which has not been
considered before. As with the arguments in
\cite{Donkin:SchurIV,Cox:blocks,JM:schaper}, the strategy is to reduce to blocks
which contain a unique maximal partition. Following Donkin, the arguments of
\cite{Donkin:SchurIV,Cox:blocks} use translation functors to make these
reductions whereas we achieve them, slightly more quickly and in greater
generality, using just the combinatorics of the Jantzen coefficients. The key
point from our perspective is to understand the partitions which contain only
horizontal hooks (\autoref{D:horizontal}). It turns out that these partitions
classify the projective indecomposable Weyl modules for the truncated Schur
algebras in characteristic zero.

To state our main result recall that a \textbf{partition} $\mu$ of $r$ is a
non-increasing sequence of non-negative integers such that
$|\mu|=\mu_1+\mu_2+\dots=r$. If $\mu$ is a partition then the
\textbf{length} of~$\mu$ is the smallest integer $\ell(\mu)$ such
that $\mu_i=0$ for $i>\ell(\mu)$. If $\ell=\ell(\mu)$ then we omit trailing zeros
and write $\mu=(\mu_1,\dots,\mu_\ell)$.

If $\lambda$ and $\mu$ are two partitions of $r$ then $\lambda$
\textbf{dominates} $\mu$, and we write $\lambda\gedom\mu$, if
$$\sum_{i=1}^s\lambda_i\ge\sum_{i=1}^s\mu_i,$$ for all $s\ge0$.
Let $\Lambda_r$ be the set of all partitions of $r$.
A subset $\Lambda$ of $\Lambda_r$ is \textbf{cosaturated}
if whenever $\lambda\in\Lambda$ and $\mu\gedom\lambda$ then
$\mu\in\Lambda$.

Let $k$ be a field of characteristic $p\ge0$ and suppose that $0\ne q\in k$ and that
$\Lambda$ is a cosaturated set of partitions of $r$.  For each partition $\mu$ of
$r$ there is a \textit{permutation module} $M(\mu)$ for the Iwahori-Hecke algebra
$\H_{k,q}(r)$ of the symmetric group $\Sym_r$. (For more details see, for example,
\cite[Chapter~3]{M:ULect}.)  The \textbf{truncated $q$-Schur algebra} with parameter
$q\in k$ and weight poset $\Lambda$ is the endomorphism algebra
$$\S_{k,q}(\Lambda)
=\End_{\H_{k,q}(r)}\(\bigoplus_{\mu\in\Lambda}M(\mu)\).$$
The algebra $\S_{k,1}(\Lambda_r)$ is Morita equivalent to the  ``classical''
Schur algebra~\cite{Green}.

The Schur algebras $\S_{k,q}(\Lambda)$ are cellular and quasi-hereditary by
\cite[Exercise~4.13]{M:ULect}, with weight poset $\Lambda$ ordered by dominance.
Thus, for each partition $\mu\in\Lambda$ there is a Weyl module, or standard
module, $\Delta_k^\mu=\Delta_k^\mu(\Lambda)$. When $\S_{k,q}(\Lambda)$ is
semisimple the Weyl modules $\set{\Delta_k^\mu|\mu\in\Lambda}$ are a complete
set of pairwise non-isomorphic irreducible $\S_{k,q}(\Lambda)$-modules.

To describe the blocks of
$\S_{k,q}(\Lambda)$ let $e\in\{0\}\cup\{2,3,\dots\}$ be minimal such that
$1+q+\dots+q^{e-1}=0$ and set $e=0$ if no such integer exists. Set
$\P=\{1,e,ep,ep^2,\dots\}$. If $\mu\in\Lambda$ let $\kappa=\core(\mu)$ be its $e$-core
(see Section~3.2). Define a \textbf{length function}
$\ell_\Lambda\map\Lambda\N$ by setting
$$\ell_\Lambda(\mu)=\min\set{i\ge0|\lambda_j=\kappa_j \text{ whenever }j>i
\text{ and } \lambda\in\Lambda \text{ has $e$-core }\kappa},$$
where $\kappa=(\kappa_1,\kappa_2,\dots)$.
If $\ell_\Lambda(\mu)\le1$ set $s_\Lambda(\mu)=1$ and otherwise define
$$s_\Lambda(\mu)=\max\set{s\in\P| \mu_i-\mu_{i+1}\equiv-1
                \pmod s \text{ for }1\le i<\ell_\Lambda(\mu)}.
$$
Finally, let $s=s_\Lambda(\mu)$ and define
$$\chi_\Lambda(\mu)=\big((\mu_1-\kappa_1)/s, (\mu_2-\kappa_2)/s,\dots\big).$$
It follows from \autoref{chi} below that $\chi_\Lambda(\mu)$ is a partition.

Our main result is the following.

\begin{MainTheorem}
Suppose that $\Lambda$ is a cosaturated set of partitions of $r$ and that
$\lambda,\mu\in\Lambda$. Then~$\Delta_k^\lambda$ and $\Delta_k^\mu$ are in the
same block as $\S_{k,q}(\Lambda)$-modules if and only if the following
three conditions are satisfied:
\begin{enumerate}
  \item $\lambda$ and $\mu$ have the same $e$-core;
  \item $s_\Lambda(\lambda)=s_\Lambda(\mu)$; and,
  \item if $s_\Lambda(\mu)>1$ then $\chi_\Lambda(\lambda)$ and $\chi_\Lambda(\mu)$
  have the same $p$-core.
\end{enumerate}
\end{MainTheorem}

Note that the $0$-core of the partition $\chi$ is $\chi$.

Let $\Lambda_{n,r}$ be the set of partitions of $r$ of length at most $n$,
so that $\Lambda_r=\Lambda_{r,r}$. Set $\S_{k,q}(n,r)=\S_{k,q}(\Lambda_{n,r})$. Then
Mod-$\S_{k,1}(n,r)$ is the category of homogeneous polynomial
representations of the general linear group $\text{GL}_n(k)$ of homogeneous
degree~$r$. The blocks of $\S_{k,q}(\Lambda_r)$ were classified by James and
Mathas~\cite[Theorem~4.24]{JM:schaper}. The blocks of $S_{k,q}(n,r)$ were
classified by Donkin~\cite[\Sect4]{Donkin:SchurIV} (for $q=1$), and Cox
\cite[Theorem 5.3]{Cox:blocks} (for $q\ne1$). We recover all of these results as
special cases of our Main Theorem.

Finally, we remark that this paper grew out of our attempts to understand
the blocks of the \textit{baby Hecke algebras}
$\H_\mu=\End_{\H_{k,q}(r)}(M(\mu))$, for $\mu$ a partition of $r$.
Let $\Lambda_\mu$ be the set of partitions of $r$ which dominate~$\mu$. Then
$\Lambda_\mu$ is cosaturated and
$\H_\mu\cong\varphi_\mu\S_{k,q}(\Lambda_\mu)\varphi_\mu$, where
$\varphi_\mu$ is the identity map on $M(\mu)$. Hence, there is a natural
\textit{Schur functor}
$$\textsf{F}_\mu\map{\Mod\S_{k,q}(\Lambda_\mu)}\Mod\H_\Lambda;
     X\mapsto X\varphi_\mu.$$
Our Main Theorem classifies the blocks of $\S_{k,q}(\Lambda_\mu)$, so
this gives a necessary condition for two $\H_\mu$-modules to belong to the
same block. Unfortunately, $\End_{\S_{k,q}(\Lambda_\mu)}(M(\mu))$ can be
larger than $\H_\mu$, so the image of a block of $\S_{k,q}(\Lambda_\mu)$
under $\textsf{F}_\mu$ need not be indecomposable. Consequently, we are not
able to describe the blocks of the algebras $\H_\mu$ completely.

\section{Jantzen equivalence and blocks}

This section develops a general theory for classifying blocks of
quasi-hereditary (cellular) algebras using Jantzen filtrations. This theory is
new in the sense that it does not appear in the literature, although everything
that we do is implicit in \cite{LM:AKblocks} which develops these results in
the special case of the cyclotomic Schur algebras.

We remark that it has long been known to people working in algebraic groups that
Jantzen filtrations could be used to determine the blocks, however, the fact
that the non-zero coefficients in the Jantzen sum formula actually classify the
blocks (\autoref{blocks}) surprised even experts in this field. For these
reasons we think it is worthwhile to give a self contained treatment of this
theory of quasi-hereditary (cellular) algebras.

\subsection{Cellular algebras} We start by recalling Graham and Lehrer's
definition of a cellular algebra~\cite{GL}. Fix an integral domain~$\O$.

\begin{Definition}[Graham and Lehrer~\cite{GL}]\label{cellular}
A \textsf{cell datum} for an associative $\O$-algebra $\S$ is a triple $(\Lambda,T,C)$
where $\Lambda=(\Lambda,>)$ is a finite poset,
$T(\lambda)$ is a finite set for $\lambda\in\Lambda$, and
$$C\map{\coprod_{\lambda\in\Lambda}T(\lambda)\times T(\lambda)}
                \S; (s,t)\mapsto C^\lambda_{st}$$
is an injective map (of sets) such that:
\begin{enumerate}
\item $\set{C^\lambda_{st}|\lambda\in\Lambda, s,t\in T(\lambda)}$
is an $\O$ basis of $\S$;
\item For any $x\in \S$ and $t\in T(\lambda)$ there exist
  scalars $r_{tv}(x)\in\O$ such that, for any $s\in T(\lambda)$,
$$C^\lambda_{st}x\equiv
\sum_{v\in T(\lambda)} r_{tv}(x)C^\lambda_{sv} \pmod{\S^\lambda},$$
where $\S^\lambda$ is the $\O$--submodule of $\S$ with basis
$\set{C^\mu_{yz}|\mu>\lambda\text{ and }y,z\in T(\mu)}$.
\item The $\O$--linear map determined by $*\map \S\S;
C^\lambda_{st}=C^\lambda_{ts}$, for all $\lambda\in\Lambda$ and
$s,t\in T(\lambda)$, is an anti--isomorphism of $\S$.
\end{enumerate}
Then $\S$ is a \textbf{cellular algebra} with \textbf{cellular basis}
$\set{C^\lambda_{st}|\lambda\in\Lambda \text{ and } s,t\in T(\lambda)}$.
\end{Definition}

Suppose that $(\Lambda,T,C)$ is a cell datum for an $\O$-algebra~$\S$.
Following Graham and Lehrer~\cite[\Sect2]{GL}, for each $\lambda\in\Lambda$
define the \textsf{cell module}, or \textbf{standard module},
$\Delta_\O^\lambda$ to be the free $\O$--module with basis
$\set{C^\lambda_t|t\in T(\lambda)}$ and with $\S$--action given by
$$C^\lambda_t x=
       \sum_{v\in T(\lambda)} r_{tv}(x)C^\lambda_{v},$$
where $r_{tv}(x)$ is the scalar from
\autoref{cellular}(b). As $r_{tv}(x)$ is independent of~$s$
this gives a well--defined $\S$--module structure on $\Delta_\O^\lambda$. The map
$\<\ ,\ \>_\lambda\map{\Delta_\O^\lambda\times \Delta_\O^\lambda}\O$
which is determined by
\begin{equation}\label{inner product}
\<C^\lambda_t,C^\lambda_u\>_\lambda C^\lambda_{sv}
     \equiv C^\lambda_{st}C^\lambda_{uv}\pmod{\S^\lambda},
\end{equation}
for $s, t, u, v\in T(\lambda)$, defines a symmetric bilinear form
on $\Delta_\O^\lambda$. This form is associative in the sense that
$\<ax,b\>_\lambda=\<a,bx^*\>_\lambda$, for all $a,b\in \Delta_\O^\lambda$ and
all $x\in \S$.

It follows easily from the definitions that the framework above is
compatible with base change. That is, if $A$ is a commutative $\O$-algebra
then $\set{C^\lambda_{st}\otimes 1_A|\lambda\in\Lambda, s,t\in T(\lambda)}$
is a cellular basis of the $A$-algebra $\S_A=\S\otimes_\O A$. Moreover,
$\Delta^\lambda_A\cong\Delta^\lambda_\O\otimes_\O A$ for all
$\lambda\in\Lambda$.

\subsection{Jantzen filtrations of cell modules}
In order to define the Jantzen filtrations of the standard modules we now
assume that~$\O$ is a discrete valuation ring with maximal ideal~$\p$ and
we let $K$ be the field of fractions of $\O$ and $k=\O/\p$ be the
residue field of $\O$. As remarked in the last paragraph,
$\S_K=\S\otimes_\O K$ and $\S_{k}=\S\otimes_\O k$ are cellular algebras with,
in essence, the same cell datum. In particular,
$\Delta_K^\lambda\cong \Delta_\O^\lambda\otimes_\O K$ and
$\Delta_k^\lambda\cong \Delta_\O^\lambda\otimes_\O k$, for
$\lambda\in\Lambda$.

\textit{Henceforth, we assume that $\S_K$ is a semisimple algebra}.
Equivalently, by \cite[Theorem~3.8]{GL}, we assume that the bilinear form
$\<\ ,\ \>_\lambda$ for $\Delta^\lambda_K$ is non-degenerate, for all
$\lambda\in\Lambda$. Thus, $\S_K$ is semisimple if and only if
$\Delta_K^\lambda$ is irreducible for all $\lambda\in\Lambda$. Hence, $(K,\O,k)$
is a modular system for $\S_{k}$.

For $\lambda\in\Lambda$ and $i\ge0$ define
$$J_i(\Delta_\O^\lambda)
     =\set{x\in \Delta_\O^\lambda|\<x,y\>_\lambda\in\p^i\text{ for all }
     y\in \Delta_\O^\lambda}.$$
Then, as the form $\<\ ,\ \>_\lambda$ is associative,
$\Delta_\O^\lambda=J_0(\Delta_\O^\lambda)\supseteq
     J_1(\Delta_\O^\lambda)\supseteq\dots$
is an $\S$-module filtration of $\Delta_\O^\lambda$.

\begin{Definition}
  Suppose that $\lambda\in\Lambda$. The \textbf{Jantzen filtration} of
  $\Delta_k^\lambda$ is the filtration
   $$\Delta_k^\lambda=J_0(\Delta_k^\lambda)\supseteq J_1(\Delta_k^\lambda)
        \supseteq\dots,$$
   where
   $J_i(\Delta_k^\lambda)=\(J_i(\Delta_\O^\lambda)+\p\Delta_\O^\lambda\)/\p\Delta_\O^\lambda$
   for $i\ge0$.
\end{Definition}

Notice that $J_i(\Delta_k^\lambda)=0$ for $i\gg0$ since $\Delta_k^\lambda$ is
finite dimensional.

For each $\lambda\in\Lambda$ set
$L_k^\lambda=\Delta_k^\lambda/J_1(\Delta_k^\lambda)$. By the general theory of
cellular algebras~\cite[Theorem 3.4]{GL}, $L_k^\lambda$ is either zero or
absolutely irreducible. Moreover, all of the irreducible $\S_{k}$-modules
arise uniquely in this way. (Note that $L_K^\lambda=\Delta_K^\lambda$, for
$\lambda\in\Lambda$, since $\S_K$ is semisimple.)

The definition of the Jantzen filtration makes sense for the standard
modules of arbitrary cellular algebras, however, for the next Lemma we need
to \textit{assume that $\S_{k}$ is quasi-hereditary}. By Remark~3.10
of \cite{GL}, $\S_{k}$ is quasi-hereditary if and only if
$J_1(\Delta_\O^\lambda)\ne\Delta_\O^\lambda$ for all
$\lambda\in\Lambda$.  Equivalently, $\S_{k}$ is quasi-hereditary if and only
if $L_k^\lambda\ne0$, for all $\lambda\in\Lambda$.

A subset $\Gamma$ of $\Lambda$ is \textbf{cosaturated} if $\lambda\in\Gamma$ whenever
$\lambda\in\Lambda$ and $\lambda>\gamma$ for some $\gamma\in\Gamma$. Let $\S^\Gamma$
be the subspace of $\S$ spanned by the elements
$\set{C^{\lambda}_{st}|\lambda\in\Gamma\text{ and }s,t\in T(\lambda)}$. For future
reference we note the following fact which follows easily from
\autoref{cellular} and the last paragraph.

\begin{Lemma}\label{L:cosaturated}
  Suppose that $\Gamma$ is a cosaturated subset of $\Lambda$. Then:
  \begin{enumerate}
    \item The algebra $\S/\S^\Gamma$ is a cellular algebra with cellular basis
      $$\set{C^{\lambda}_{st}+\S^\Gamma|\lambda\in\Lambda{\setminus}\Gamma
                    \text{ and }s,t\in T(\lambda)}.$$
    \item If $\S$ is a quasi-hereditary algebra then so is $\S/\S^\Gamma$.
  \end{enumerate}
\end{Lemma}

Let $K_0(\S_{k})$ be the Grothendieck group of finite dimensional right
$\S_{k}$-modules. If $M$ is an $\S_{k}$-module let $[M]$ be its image in
$K_0(\S_{k})$.

If $M$ is an $\S_{k}$-module and $\mu\in\Lambda$ let $[M:L_k^\mu]$ be the
multiplicity of the simple module $L_k^\mu$ as a composition factor of $M$.  In
particular, if $\lambda,\mu\in\Lambda$ let
$d_{\lambda\mu}=[\Delta_k^\lambda:L_k^\mu]$. Then, by
\cite[Proposition~3.6]{GL}, $d_{\lambda\lambda}=1$ and $d_{\lambda\mu}\ne0$ only
if $\lambda\ge\mu$. Consequently, the decomposition matrix
$\big(d_{\lambda\mu}\big)_{\lambda,\mu\in\Lambda}$ of $\S_{k}$ is a square
unitriangular matrix, when its rows and columns are ordered in a way that is
compatible with~$>$. Therefore, the decomposition matrix of $\S_{k}$ is
invertible over~$\Z$ and as a consequence we obtain the following.

\begin{Lemma}\label{vanishing}
  Suppose that $\S_{k}$ is a quasi-hereditary cellular algebra and
  $\lambda,\mu\in\Lambda$. Then
  \begin{enumerate}
    \item $\set{[\Delta_k^\lambda]|\lambda\in\Lambda}$ is a $\Z$-basis of
      $K_0(\S_{k})$.
    \item There exist integers
      $J_{\lambda\mu}\in\Z$ such that
      $$\sum_{i>0} [J_i(\Delta_k^\lambda)]
           =\sum_{\substack{\mu\in\Lambda\\\lambda>\mu}}
                   J_{\lambda\mu}[\Delta_k^\mu].$$
    \item If $\mu\ne\lambda$
      then $d_{\lambda\mu}=[J_1(\Delta_k^\lambda):L_k^\mu]$. Consequently,
      $d_{\lambda\mu}\ne0$ if and only
      if $J_{\lambda\mu}'\ne0$, where
      $$J_{\lambda\mu}'=\sum_{i>0}[J_i(\Delta^\lambda_k):L_k^\mu]
                       =\sum_{\lambda>\nu\ge\mu}J_{\lambda\nu}d_{\nu\mu}.$$
     Moreover, if $\mu\ne\lambda$ then $d_{\lambda\mu}\le J'_{\lambda\mu}$.
  \end{enumerate}
\end{Lemma}

The integers $J_{\lambda\mu}$ are the \textbf{Jantzen coefficients} of
$\S_{k}$. By definition,
\begin{equation}\label{J'}
J_{\lambda\mu}'=[\bigoplus_{i>0}J_i(\Delta_k^\lambda):L_k^\mu]
            \ge[\rad\Delta_k^\lambda:L_k^\mu],
\end{equation}
where $\rad\Delta_k^\lambda=J_1(\Delta_k^\lambda)$ is the radical of $\Delta_k^\lambda$. We show in
the next section that the Jantzen coefficients determine the blocks of~$\S_{k}$. They also determine the
irreducible standard modules.

\begin{Corollary}\label{simple}
  Suppose that $\lambda\in\Lambda$. Then the following are equivalent:
  \begin{enumerate}
    \item $\Delta^\lambda_k=L^\lambda_k$ is an irreducible
      $\S_{k}$-module.
    \item $d_{\lambda\mu}=\delta_{\lambda\mu}$, for all $\mu\in\Lambda$ (Kronecker delta).
    \item $J_{\lambda\mu}=0$, for all $\mu\in\Lambda$.
  \end{enumerate}
\end{Corollary}

\begin{proof}Parts (a) and (b) are equivalent essentially by definition and
  (b) and (c) are equivalent by \autoref{vanishing}(c).
\end{proof}

\subsection{Jantzen coefficients and the blocks of $\S_{k}$}
The algebra $\S_{k}$ decomposes in a unique way as a direct sum of
indecomposable two-sided ideals $\S_{k}=B_1\oplus\dots\oplus B_d$. If $M$ is an
$\S_{k}$-module then $MB_i$ is a $B_i$-module. We say that $M$
\textbf{belongs} to the block $B_i$ if $MB_i=M$. Using an idempotent argument
(cf.~\cite[Theorem~56.12]{C&R:2}) it is easy to show that two indecomposable
$\S_{k}$-modules $P$ and~$Q$ belong to the same block if and only if they are
in the same \emph{linkage class}. That is, there exist indecomposable
modules $P_1=P,\dots,P_l=Q$ such that $P_i$ and~$P_{i+1}$ have a common
irreducible composition factor, for $i=1,\dots,l-1$.

\begin{Definition}
  Suppose that $\lambda,\mu\in\Lambda$. Then $\lambda$ and $\mu$ are
  \textbf{Jantzen equivalent}, and we write $\lambda\sim_J\mu$, if there
  exist $\lambda_1=\lambda,\lambda_2,\dots,\lambda_l=\mu\in\Lambda$ such
  that either
  $$J_{\lambda_i\lambda_{i+1}}\ne0\quad\text{or}\quad
    J_{\lambda_{i+1}\lambda_i}\ne0,$$
  for $1\le i<l$.
\end{Definition}

The next result shows that the Jantzen equivalence classes and the blocks of $\S_{k}$ coincide. This is
the main result of this section.

\begin{Proposition}\label{blocks}
  Suppose that $\lambda,\mu\in\Lambda$. Then $\Delta_k^\lambda$ and
  $\Delta_k^\mu$ belong to the same block as $\S_{k}$--modules if
  and only if~$\lambda\sim_J\mu$.
\end{Proposition}

\begin{proof}
We essentially repeat the argument of \cite[Proposition~2.9]{LM:AKblocks}.
Before we begin observe that if $\nu\in\Lambda$ then $\Delta_k^\nu$ is
indecomposable because $L_k^\nu$ is the simple head of $\Delta_k^\nu$.
Consequently, all of the composition factors of $\Delta_k^\nu$ belong to the
same block.

Suppose, first, that $\lambda\sim_J\mu$. By definition
$J_i(\Delta_k^\lambda)$ is a submodule of $\Delta_k^\lambda$ for all~$i$,
so all of the composition factors of $\bigoplus_{i>0}J_i(\Delta_k^\lambda)$
belong to the same block as $\Delta_k^\lambda$ by the last paragraph.
Define $\Lambda'$ to be the subset of $\Lambda$ such that $\nu\in\Lambda'$
whenever $\Delta_k^\nu$ and $\Delta_k^\lambda$ are in different blocks.
Then $\sum_{\nu\in\Lambda'}J_{\lambda\nu}[\Delta_k^\nu]=0$ by \autoref{J'}.
Hence, $J_{\lambda\nu}=0$ whenever $\nu\in\Lambda'$ by \autoref{vanishing}(a).
It follows that $\Delta_k^\lambda$ and $\Delta_k^\mu$ belong to the same
block whenever $J_{\lambda\mu}\ne0$.

To prove the converse it is sufficient to show that $\lambda\sim_J\mu$
whenever $d_{\lambda\mu}\ne0$. By \autoref{vanishing}(c), if
$d_{\lambda\mu}\ne0$ then $J_{\lambda\mu}'\ne0$. Therefore, there exists a
partition $\nu_1\in\Lambda$ such that $\lambda>\nu_1\ge\mu$,
$J_{\lambda\nu_1}\ne0$ and $d_{\nu_1\mu}\ne0$. If $\nu_1=\mu$ then
$\lambda\sim_J\mu$ and we are done. If $\nu_1\ne\mu$ then
$d_{\nu_1\mu}\ne0$, so $J_{\nu_1\mu}'\ne0$ and we may repeat this argument
to find $\nu_2\in\Lambda$ with $\nu_1>\nu_2\ge\mu$, $J_{\nu_1\nu_2}\ne0$
and $d_{\nu_2\mu}\ne0$.  Continuing in this way we can find elements
$\nu_0=\lambda,\nu_1,\dots,\nu_l=\mu$ in $\Lambda$ such that
$J_{\nu_{i-1}\nu_i}\ne0$, $d_{\nu_i\mu}\ne0$, for $0< i<l$, and
$\lambda>\nu_1>\dots>\nu_l=\mu$. Note that we must have
$\nu_l=\mu$ for some $l$ since $\Lambda$ is finite. Therefore,
$\lambda\sim_J\nu_1\sim_J\dots\sim_J\nu_l=\mu$ as required.
\end{proof}

We have chosen to prove \autoref{blocks} using the formalism of
cellular algebras, however, it can be proved entirely within the framework
of quasi-hereditary algebras. Suppose that $\O$ is a (complete) discrete
valuation ring with residue field~$k$. Let $A$ be a quasi-hereditary algebra
which is free and of finite rank as an $\O$-module and set
$A_k=A\otimes_\O k$. Following, for example,
McNinch~\cite[\Sect4.1]{McNinch:Howe} we can define Jantzen filtrations of
the standard modules of the quasi-hereditary algebra $A_k$. The standard
modules of a quasi-hereditary algebra are always indecomposable and they
always give a basis for the Grothendieck group of~$A_k$. Moreover, the
decomposition matrix of~$A_k$ is unitriangular. Using these general facts,
\autoref{blocks} (and \autoref{vanishing} and
\autoref{simple}), can be proved for $A_k$ following the arguments
above.

\section{Combinatorics and Jantzen equivalence for Schur algebras}

We are now ready to start proving our Main Theorem. We begin by recalling the
combinatorics we need to describe the Jantzen coefficients for the algebras
$\S_{k,q}(\Lambda)$ from the introduction. Recall from the introduction that
$\Lambda_r=\Lambda_{r,r}$ is the set of all partitions of~$r$.

As in the introduction, let $\Lambda$ be a cosaturated set of partitions of $r$
and fix a field $k$ of characteristic $p\ge0$ and a non-zero element $q\in
k^\times$. Let $\S_{k,q}(\Lambda)$ be the $q$-Schur algebra over $k$ with
parameter $q$ and weight poset $\Lambda$.

\subsection{Schur functors}
Recall that $\Lambda_r$ is the set of all partitions of $r$ and that
$\S_{k,q}(\Lambda_r)$ is the $q$-Schur algebra with weight poset
$\Lambda_r$. There is a natural embedding
$\S_{k,q}(\Lambda)\hookrightarrow\S_{k,q}(\Lambda_r)$. Moreover, it is easy to see that if
$e_\Lambda$ is the identity element of $\S_{k,q}(\Lambda)$ then
$\S_{k,q}(\Lambda)=e_\Lambda\S_{k,q}(\Lambda_r)e_\Lambda$.

For the next result, write $\Delta_k^\lambda(\Lambda_r)$ for the standard
modules of $\S_{k,q}(\Lambda_r)$ and $\Delta_k^\lambda(\Lambda)$ for the
standard modules of $\S_{k,q}(\Lambda)$. Then, by standard arguments (see for
example, \cite[\Sect6]{Green} or \cite[Proposition A3.11]{Donkin:book}), we obtain the
following.

\begin{Lemma}[Schur functor]\label{Schur functor}
  Right multiplication by $e_\Lambda$ induces an exact functor
  $F_\Lambda$ from the category of right
  $\S_{k,q}(\Lambda_r)$--modules to the category of right $\S_{k,q}(\Lambda)$-modules
  such that
  $$F_\Lambda(\Delta^\lambda_k(\Lambda_r))\cong\begin{cases}
     \Delta^\lambda_k(\Lambda),&\text{if $\lambda\in\Lambda$},\\
     0,&\text{if $\lambda\notin\Lambda$}.
   \end{cases}$$
  Moreover, if $\lambda\in\Lambda$ then
  $\dim \Delta^\lambda_k(\Lambda)=\dim\Delta^\lambda_k(\Lambda_r)$.
\end{Lemma}

The standard modules $\Delta_k^\lambda(\Lambda)$ are often called the
\textbf{Weyl modules} of $\S_{k,q}(\Lambda)$. In view of \autoref{Schur
functor} we now write $\Delta^\lambda_k=\Delta^\lambda_k(\Lambda)$.

Define the \textbf{quantum characteristic} of $(k,q)$ to be
$$e=\min\set{c\ge1|1+q+\dots+q^{c-1}=0}$$
and set $e=0$ if no such integer exists. That is, $e=p$ if $q=1$, $e=0$ if $q$
is not a root of unity and otherwise $e$ is the multiplicative order of~$q$.

By \cite[Exercise~4.14]{M:ULect} the algebra $\S_{k,q}(\Lambda_r)$ is semisimple if
and only if $e=0$ or $e>r$. Hence, applying the Schur functor of
\autoref{Schur functor}, $\S_{k,q}(\Lambda)$ is semisimple if $e=0$ or $e>r$.
Consequently, we assume henceforth that $0<e\le r$.

\subsection{Jantzen  coefficients for $\S_{k,q}(\Lambda)$}
To define the Jantzen filtrations of $\S_{k,q}(\Lambda)$ fix a modular system
\textit{with parameters} $(K,\O,k)_{t,q}$ such that
\begin{itemize}
  \item $\O$ is a discrete valuation ring with maximal ideal $\p$ and
    $t$ is an invertible element of~$\O$;
  \item $K$ is the field of fractions of $\O$ and
$\S_{K,t}(\Lambda)\cong S_{\O,t}(\Lambda)\otimes_\O K$ is semisimple; and,
  \item $k\cong\O/\p$, $q=t+\p\O$ and
    $\S_{k,q}(\Lambda)\cong \S_{\O,t}(\Lambda)\otimes_\O k$.
\end{itemize}
In general, the Jantzen filtrations of $\S_{k,q}(\Lambda)$-modules may depend
upon this choice of modular system. In this paper, however, we only need to know
whether or not the Jantzen coefficients are zero and this is independent of the
choice of modular system by \autoref{non-vanishing rows} below.

Least the reader be concerned that a modular system with these properties need
not always exist we note that if $x$ is an indeterminate over $k$ then we could
let $\O=k[x]_{(x)}$ be the localization of $k[x]$ at the prime ideal $(x)$, so
that $\p=x\O$ is the unique maximal ideal of $\O$. Then $\O$ is a discrete
valuation ring with field of fractions $K=k(x)$. Set $t=x+q$ which, by abuse of
notation, we consider as an invertible element of~$\O$. Then, using the remarks
above, the reader can check that $(K,\O,k)_{t,q}$ is a modular system with
parameters for $\S_{k,q}(\Lambda)$.

As in section~2, for each partition $\lambda\in\Lambda$ define the Jantzen
filtration $\{J_i(\Delta_k^\lambda)\}$ of the standard module
$\Delta_k^\lambda$ of~$\S_{k,q}(\Lambda)$. Define the
\textbf{Jantzen coefficients} of $\S_{k,q}(\Lambda)$ to be the integers
$J_{\lambda\mu}^\Lambda$ determined by the following equations in
$K_0(\S_{k,q}(\Lambda))$:
$$\sum_{i>0}[J_i(\Delta^\lambda_k)]
    =\sum_{\substack{\mu\in\Lambda\\\mu>\lambda}}
        J_{\lambda\mu}^\Lambda[\Delta^\mu_k].$$
Recall that $\Lambda_r$ is the
set of partitions of $r$. For $\lambda,\mu\in\Lambda_r$, set
$J_{\lambda\mu}=J_{\lambda\mu}^{\Lambda_r}$. Applying the Schur functor
(\autoref{Schur functor}), shows that the Jantzen coefficients depend
only on $\Lambda_r$ in the following sense.

\begin{Corollary}
  Suppose that $\lambda,\mu\in\Lambda$. Then
  $J_{\lambda\mu}^\Lambda=J_{\lambda\mu}$.
\end{Corollary}

Henceforth, we write $J_{\lambda\mu}^\Lambda=J_{\lambda\mu}$ for the
Jantzen coefficients of $\S_{k,q}(\Lambda)$.

\subsection{Beta numbers, abaci and cores}
We now introduce the notation that we need to describe when the Jantzen
coefficients of $\S_{k,q}(\Lambda)$ are non-zero. The bulk of the work has already been
done in~\cite{JM:schaper}.

For any partition $\mu=(\mu_1,\mu_2,\dots)$ let
$\Diag(\mu)=\set{(i,j)|1\le j\le\mu_i}$ be the \textbf{diagram} of~$\mu$ which
we think of as a (left justified) collection of boxes in the plane. The
\textbf{$e$-residue} of a node $(i,j)\in\Diag(\mu)$ is the unique integer
$\res_e(i,j)$ such that $0\le\res_e(i,j)<e$ and $\res_e(i,j)\equiv j-i\pmod e$.

Fix any integer $l\ge\ell(\mu)$. For $1\le i\le l$ set
$\beta_i=\mu_i-i+l$. Then $\beta_1>\beta_2>\dots>\beta_l\ge0$ are the
\textbf{$l$-beta numbers} for $\mu$. It is well-known and easy to prove
that the beta numbers give a bijection between the set of partitions with
at most $l$ non-zero parts and the set of strictly increasing non-negative
integer sequences of length~$l$.

An \textbf{$e$-abacus}~\cite{James:YoungD} is a Chinese abacus with $e$
\textbf{runners} and with \textbf{bead positions} numbered $0,1,2,\dots$ from left to
right and then top to bottom.  (We will also need $p$-abaci and $s$-abaci.) Let
$\beta_1>\beta_2>\dots>\beta_l$ be the sequence of $l$-beta numbers for $\lambda$.
The $l$-bead \textbf{abacus configuration} for $\lambda$ is the abacus with
beads at positions $\beta_1,\beta_2,\dots,$ and $\beta_l$.  Any abacus
configuration determines a set of beta numbers and hence corresponds to a unique
partition. If $\beta\ge0$ then bead position $\beta+1$ is the position to the
\textbf{right} of $\beta$ and bead position $\beta-1$ (if $\beta>0$) is the
position to the \textbf{left}. (In particular, the bead position to the
left of a position on runner $0$ is on runner $e-1$ in the previous row.)

For example, taking $e=3$ and $l=6$ the abacus configurations for the
partitions $\lambda=(4,4,3,1)$ and $\kappa=(4,2)$ are as follows:
$$\abacus35{0,1,3,6,8,9}\hspace*{20mm}
  \abacus35{0,1,2,3,6,9}$$

The \textbf{$e$-core} of $\lambda$ is the partition $\core(\lambda)$ whose abacus
configuration is obtained from an abacus configuration of~$\lambda$ by pushing all
beads as high as possible on their runner. If $e=0$ then, by convention,
$\core[0](\lambda)=\lambda$. If $e>0$ then the \textbf{$e$-weight} of $\lambda$ is
the integer $(|\lambda|-|\core(\lambda)|)/e$ otherwise $\lambda$ has $e$-weight zero.
For example, if $\lambda=(4,4,3,1)$, as in the example above, then
$\core(\lambda)=(4,2)=\kappa$ and~$\lambda$ has $3$-weight~$2$.

Observe that, up to a constant shift, $\set{\lambda_i-i|1\le i\le l}$ is
the set of beta numbers of~$\lambda$. Therefore, two partitions $\lambda$
and $\mu$ of~$r$ have the same $e$-core if and only if
$$\lambda_i-i\equiv\mu_{i^w}-i^w\pmod e,$$
for some $w\in\Sym_r$.

Let $\lambda'=(\lambda_1',\lambda_2',\dots)$ be the partition
\textbf{conjugate} to $\lambda$, so that
$\lambda_j'=\#\set{i\ge1|\lambda_i\ge j}$ for $j\ge1$. If
$(a,b)\in\Diag(\lambda)$ then the \textbf{$(a,b)$-rim hook} of~$\lambda$ is the
set of nodes
$$R^\lambda_{ab}=\set{(i,j)\in\Diag(\lambda)|
         a\le i\le\lambda_b', b\le j\le\lambda_a \text{ such that }
         (i+1,j+1)\notin\Diag(\lambda)}.$$
The node $(a,\lambda_a)$ is the \textbf{hand node} of $R^\lambda_{ab}$,
$f^\lambda_{ab}=(\lambda_b',b)$ is the \textbf{foot node} and
$h^\lambda_{ab}=|R^\lambda_{ab}|$ is the \textbf{hook length} of
$R^\lambda_{ab}$. The \textbf{foot residue} of $R^\lambda_{ab}$ is
$b-\lambda_b'\pmod e$, the residue of $f^\lambda_{ab}$.  The hook
$R^\lambda_{ab}$ is an \textbf{$h$-hook} if $h=h^\lambda_{ab}$.  Thus,
$R^\lambda_{ab}$ is the $h$-hook consisting of the set of nodes along the
`rim' of $\Diag(\lambda)$ which connects the hand and foot nodes.

If $\mu$ is a partition and $\Diag(\mu)=\Diag(\lambda)\setminus R^\lambda_{ab}$, for some $(a,b)\in\Diag(\lambda)$, then
we say that $\mu$ is obtained from $\lambda$ by \textbf{unwrapping} the rim hook $R^\lambda_{ab}$ and that
$\lambda$ is obtained from $\mu$ by \textbf{wrapping on} this hook. A hook $R^\lambda_{ab}$ is
\textbf{removable} if $\Diag(\lambda)\setminus R^\lambda_{ab}$ is the diagram of a partition.

Using the definitions, if $\{\beta_1,\dots,\beta_l\}$ is a set of beta
numbers for $\lambda$ and if $\mu$ is obtained from $\lambda$ by unwrapping
the rim hook $R^\lambda_{ab}$ then
$\{\beta_1,\dots,\beta_{a-1},\beta_a-h^\lambda_{ab},\beta_{a+1},\dots,\beta_l\}$
is a set of beta numbers for $\mu$. With a little extra care we obtain the
following well-known fact; see, for example, \cite[Lemma 5.26]{M:ULect}.

\begin{Lemma}\label{wrapping}
Suppose that $\lambda$ is a partition. Then moving a bead $h$ positions to the
right in the abacus configuration of $\lambda$ from runner~$f$ to runner $f'$
corresponds to wrapping an $h$--rim hook with foot residue $f$ onto $\lambda$.
Equivalently, moving a bead~$h$ positions to the left, from runner $f'$ to runner
$f$ corresponds to unwrapping an $h$--rim hook from~$\lambda$ with foot residue $f$.
\end{Lemma}

The following Lemma shows that the non-vanishing of the Jantzen
coefficients is independent of the choice of modular system. First,
however, we need a definition.  Let $\nu_{e,p}\map{\N\setminus\{0\}}\N$ be the map which
sends $h>0$ to
$$\nu_{e,p}(h)=\begin{cases}
  \nu_p(\frac he)+1, &\text{if }e\mid h,\\
  0,&\text{if }e\nmid h,
\end{cases}$$
where $\nu_p$ is the usual $p$-adic valuation map on $\N\setminus\{0\}$ (and
where we set $\nu_p(h)=0$ for all $h\in\N\setminus\{0\}$ if $p=0$). Note that
$\nu_{p,p}(h)=\nu_p(h)$, for all $h\in\N\setminus\{0\}$.

For each integer $k\in\Z$, define the \textbf{$t$-quantum integer}
$[k]=(t^k-1)/(t-1)\in\N[t,t^{-1}]$.

The next result is a sharpening of results from \cite{JM:schaper}. The main
point is to show that the Jantzen coefficients depend only on~$e$ and~$p$.

\begin{Proposition}\label{non-vanishing rows}
  Suppose that $\lambda,\mu\in\Lambda$. Then $J_{\lambda\mu}\ne0$ if only
  if $\lambda\gdom\mu$ and there exist nodes $(a,b),(a,c)\in\Diag(\lambda)$ such
  that $b<c$ and $\nu_{e,p}(h^\lambda_{ab})\ne\nu_{e,p}(h^\lambda_{ac})$
  and $\mu$ is obtained from $\lambda$ by unwrapping the rim hook
  $R^\lambda_{ac}$ and then wrapping it back on with its hand node in
  column~$b$.
\end{Proposition}

\begin{proof}
  Let $\nu_\p$ be the $\p$-adic valuation map on~$\O$. Then,
  by \cite[Theorem 4.3]{JM:schaper},
  $$J_{\lambda\mu}=\sum_{(a,b),(a,c)\in\Diag(\lambda)}
          \pm\(\nu_\p([h^\lambda_{ab}])-\nu_\p([h^\lambda_{ac}])\),$$
  where the sum is over a collection of nodes $(a,b),(a,c)\in\Diag(\lambda)$
  which satisfy the assumptions of the Lemma. Using the abacus and
  \autoref{wrapping}, it is easy to see that there is at most one pair of
  nodes $(a,b),(a,c)\in\Diag(\lambda)$ that allow us to obtain $\mu$ from
  $\lambda$ by unwrapping and then wrapping on a hook, so the last equation
  becomes
  $$J_{\lambda\mu}=\pm\(\nu_\p([h^\lambda_{ab}])-\nu_\p([h^\lambda_{ac}])\).$$
  (The sign is determined by the parity of the sum of the leg lengths of
  the rim hooks involved.) By \cite[Lemma~4.17]{JM:schaper}, if $a$ and $b$
  are any integers then $\nu_\p([a])=\nu_\p([b])$ if and only if either (a)
  $e\nmid a$ and $e\nmid b$, or (b) $e\mid a$, $e\mid b$ and
  $\nu_p(a)=\nu_p(b)$.
  Putting these two statements together proves the lemma.
\end{proof}

Surprisingly, the next result appears to be new.

\begin{Corollary}\label{non-vanishing cols}
  Suppose that $\lambda,\mu\in\Lambda$. Then $J_{\lambda\mu}\ne0$ if and
  only if $\lambda\gdom\mu$ and there exist nodes $(x,z),(y,z)\in\Diag(\mu)$
  such that $x<y$ and
  $\nu_{e,p}(h^\mu_{xz})\ne\nu_{e,p}(h^\mu_{yz})$ and $\lambda$ is
  obtained from $\mu$ by unwrapping the rim hook $R^\mu_{yz}$ and
  then wrapping it back on with its foot node in row~$x$.
\end{Corollary}

\begin{proof}
  Fix nodes $(a,b),(a,c)\in\Diag(\lambda)$ as in \autoref{non-vanishing
  rows} such that $b<c$ and $\mu$ is obtained from~$\lambda$ by unwrapping
  $R^\lambda_{ac}$ and wrapping it back on with its hand node in column~$b$. Let
  $(y,z)\in\Diag(\mu)$ be the unique node such that
  $\Diag(\mu){\setminus}R^\mu_{yz}=\Diag(\lambda){\setminus}R^\lambda_{ac}$, as
  in the diagram below.
  \begin{figure}[h]
  $$\begin{tikzpicture}[scale=0.25,every node/.append style={font=\footnotesize}]
      \useasboundingbox (0,0) rectangle (13,-10);
      \draw[very thin](-0.5,-1.75)node[left]{$a$}--(11,-1.75);
      \draw[very thin](5.25,0.5)node[above]{$b$}--(5.25,-7);
      \draw[very thin](7.75,0.5)node[above]{$c$}--(7.75,-4);
      \fill[blue!70](11,-1.5)rectangle(8.5,-2);
      \fill[blue!70](9,-1.5)rectangle(8.5,-4);
      \fill[blue!70](9,-3.5)rectangle(7.5,-4);
      \fill[red!70](5.5,-6)rectangle(5,-7.5);
      \fill[red!70](5.5,-7)rectangle(3,-7.5);
      \fill[red!70](3.5,-7)rectangle(3,-9.5);
      \draw[very thin](-0.5,-3.75)node[left]{$x$}--(7.5,-3.75);
      \draw[very thin](-0.5,-6.25)node[left]{$y$}--(5,-6.25);
      \draw[very thin](3.25,0.5)node[above]{$z$}--(3.25,-7);
      \node(Rac)at(10.25,-3){$R^\lambda_{ac}$};
      \node(Ryz)at(4.75,-8.5){$R^\mu_{yz}$};
      \node(0,0)[above left]{\Large$\lambda$};
      \path[thick,dashed,->](Rac) edge[bend left](Ryz);
      \draw[very thick](0,0)--++(13,0)--++(0,-1)--++(-2,0)--++(0,-1)--++(-2,0)--++(0,-2)
            --++(-2,0)--++(0,-2)--++(-2,0)--++(0,-1)--++(-2,0)--++(0,-3.5)
            --++(-1.5,0)--++(0,-1)--++(-1.5,0)--cycle;
  \end{tikzpicture}$$
  \end{figure}
  Thus, $y=\lambda_b'+1$ and~$R^\mu_{yz}$ is the rim hook which is wrapped back
  on to $\Diag(\lambda)\setminus R^\lambda_{ac}$ to form~$\mu$. Set
  $x=\lambda_c'$. Then $x<y$ and
  $R^\mu_{xz}\sqcup R^\lambda_{ac}=R^\lambda_{ab}\sqcup R^\mu_{yz}$ (disjoint
  unions), where these sets are disjoint because $b<c$. Therefore,
  $h^\mu_{xz}=h^\lambda_{ab}$, since $h^\mu_{yz}=h^\lambda_{ac}$, and $\lambda$
  is obtained by unwrapping $R^\mu_{yz}$ from $\mu$ and then wrapping it back on
  with its foot node in row~$x$. The result now follows by
  \autoref{non-vanishing rows}.
\end{proof}

Hence, the coefficients $J_{\lambda\mu}$ are almost symmetric in $\lambda$ and
$\mu$.

\begin{Corollary}\label{symmetry}
  Suppose that $\lambda,\lambda',\mu,\mu'\in\Lambda_r$. Then
  $J_{\lambda\mu}\ne0$ if and only if $J_{\mu'\lambda'}\ne0$.
\end{Corollary}

\begin{proof}
  This is immediate because the conditions in \autoref{non-vanishing rows} and
  \autoref{non-vanishing cols} are interchanged by taking conjugates of partitions.
\end{proof}

It is well-known from~\cite{JM:schaper} that the blocks of~$\S_{k,q}(n,n)$ are
determined by the $e$-cores of the partitions. The next lemma establishes the
easy half of this classification within our framework.

\begin{Lemma}\label{e-cores}
  Suppose that $\lambda,\mu\in\Lambda$ and that $J_{\lambda\mu}\ne0$. Then
  $\lambda\gdom\mu$ and $\lambda$ and $\mu$ have the same $e$-core.
\end{Lemma}

\begin{proof}By \autoref{non-vanishing rows} since $J_{\lambda\mu}\ne0$ there
  exist nodes $(a,b),(a,c)\in\Diag(\lambda)$  such that $b<c$,
  $\nu_{e,p}(h^\lambda_{ab})\ne\nu_{e,p}(h^\lambda_{ac})$ and $\mu$ is
  obtained from $\lambda$ by unwrapping the rim hook $R^\lambda_{ac}$ and
  then wrapping it back on with its hand node in column~$b$. As the rim
  hook $R^\lambda_{ac}$ is wrapped back onto $\lambda$ lower down
  it follows $\lambda\gdom\mu$. It remains to show that $\lambda$ and $\mu$ have the
  same $e$-core.

  By \autoref{wrapping},  the abacus configuration for $\mu$ is obtained
  from the abacus configuration for $\lambda$ by moving one bead
  $h^\lambda_{ac}$ positions to the left and another bead $h^\lambda_{ac}$
  positions to the right.

  If $e\mid h^\lambda_{ac}$ then, by \autoref{wrapping}, the $e$-abacus
  configuration for $\mu$ is obtained by moving two beads on the same runner
  in the $e$-abacus configuration for $\mu$. Hence, $\lambda$ and $\mu$
  have the same $e$-core.

  On the other hand, if $e\nmid h^\lambda_{ac}$ then
  $\nu_{e,p}(h^\lambda_{ac})=0$, so that $e$ divides $h^\lambda_{ab}$ by
  \autoref{non-vanishing rows}. Therefore, the foot residues of the
  hooks being unwrapped and then wrapped back
  into~$\lambda{\setminus}R^\lambda_{ac}$ coincide (since, modulo~$e$, these
  residues differ by $h^\lambda_{ab}$). Hence, applying \autoref{wrapping}, the
  abacus configuration for $\mu$ is obtained from the abacus for $\lambda$ by
  moving a bead $h^\lambda_{ab}$ positions to the left to runner~$f$, say, and
  then moving another bead~$h^\lambda_{ac}$ positions to the right from
  runner~$f$. Consequently, the number of beads of each runner is unchanged, so
  that~$\lambda$ and~$\mu$ have the same $e$-core.
\end{proof}

One consequence of \autoref{e-cores} is that we can weaken the definition
of a cosaturated set of partitions.

\begin{Definition}\label{D:Cosaturated}
  Suppose that $\Lambda$ is a set of partitions of~$r$. Then $\Lambda$ is
  \textbf{$e$-cosaturated} if $\mu\in\Lambda$ whenever there exists a partition
  $\lambda\in\Lambda$ such that $\mu\gedom\lambda$ and $\lambda$ and $\mu$ have
  the same $e$-core.
\end{Definition}

Suppose that $\Lambda$ is a cosaturated set of partitions and that $\kappa$
is an $e$-core. Let $\Lambda_\kappa$ be the set of partitions in $\Lambda$
which have $e$-core $\kappa$. Then $\Lambda_\kappa$ is $e$-cosaturated.
For each $\lambda\in\Lambda$ let $\varphi_\lambda$ be the identity map on
$M(\lambda)$ and set
$$e_\kappa^{(r)}=\sum_{\lambda\in\Lambda_\kappa}\varphi_\lambda.$$
Then $e_\kappa^{(r)}$ is an idempotent in $\S_{k,q}(\Lambda)$. Hence, \autoref{e-cores}
and standard Schur functor arguments, as in~\autoref{Schur functor},
imply that the algebra
\begin{equation}\label{E:SLambdakappa}
 \S_{k,q}(\Lambda_\kappa):= e_\kappa^{(r)} \S_{k,q}(\Lambda) e_\kappa^{(r)},
\end{equation}
is a quasi-hereditary algebra with weight poset $\Lambda_\kappa$. The
algebra $\S_{k,q}(\Lambda_\kappa)$ is a direct sum of blocks of $\S_{k,q}(\Lambda)$. In
general, however, $\S_{k,q}(\Lambda_\kappa)$ is not indecomposable. In what follows it will
sometimes be convenient to assume that $\Lambda=\Lambda_\kappa$ is
$e$-cosaturated.

\subsection{Projective $\S_{k,q}(\Lambda)$-modules} Up until now we have been
recalling and slightly improving on results from the literature, but we now
leave this well-trodden path. The main result of this section is \autoref{conditions}
which is a very subtle characterisation of the partitions which contain only
\textit{horizontal} $e$-hooks. This result is the key to the main results of
this paper. In particular, it motivates the combinatorial definitions which
underpin our Main Theorem.

Let $P^\mu_k$ be the projective cover of $L_k^\mu$. In \autoref{simple} we used
the Jantzen coefficients to classify the simple standard modules of a
quasi-hereditary algebra. The next result, which is routine but puts the results
below into context, shows that the Jantzen coefficients also classify the
projective standard modules of $\S_{k,q}(\Lambda)$.

\begin{Lemma}\label{projective}
  Suppose that $\mu\in\Lambda$. Then the following are equivalent:
  \begin{enumerate}\topsep=0pt
    \item $\Delta_k^\mu=P_k^\mu$ is a projective $\S_{k,q}(\Lambda)$-module.
    \item $\Delta_k^\mu$ is a projective $\S_{k,q}(\Lambda_r)$-module.
    \item $d_{\lambda\mu}=\delta_{\lambda\mu}$, for all $\lambda\in\Lambda$.
    \item $J_{\lambda\mu}=0$, for all $\lambda\in\Lambda$.
    \item $\nu_{e,p}(h^\mu_{ac})=\nu_{e,p}(h^\mu_{bc})$, for all
      nodes $(a,c),(b,c)\in\Diag(\mu)$.
    \item $\Delta^{\mu'}_k=L_k^{\mu'}$ is an irreducible $\S_{k,q}(\Lambda_r)$-module.
  \end{enumerate}
\end{Lemma}

\begin{proof}
  Parts (a), (c) and (d) are equivalent exactly as in \autoref{simple}.
  Let $P^\mu_k$ be the projective cover of~$L^\mu_k$.

  For part~(c), it follows from the general theory of cellular
  algebras~\cite[\Sect3]{GL} that $P^\mu_k$ has a $\Delta$-filtration in which
  $\Delta^\lambda_k$ appears with multiplicity~$d_{\lambda\mu}$. (Note that if
  $d_{\lambda\mu}\ne0$ then $\lambda\gedom\mu$ so that $\lambda\in\Lambda$ since
  $\Lambda$ is cosaturated.) Consequently, $P^\mu_k$ is also the projective
  cover of $\Delta^\mu_k$, so that $\Delta^\mu_k$ is projective if and only if
  $d_{\lambda\mu}=\delta_{\lambda\mu}$. Hence, parts~(b) and~(c) are also
  equivalent.

  Finally, note that $J_{\lambda\mu}\ne0$ if and only if $J_{\mu'\lambda'}\ne0$
  by \autoref{symmetry}, so that~(d) and~(f) are equivalent by
  \autoref{simple} and (d) and (e) are equivalent by
  \autoref{non-vanishing rows}. This completes the proof.
\end{proof}

Parts (a)--(d) of \autoref{projective} are equivalent for any
quasi-hereditary cellular algebra.

\begin{Example}\label{first example}
  Suppose that $(e,p)=(3,0)$ and let $\Lambda$ be the set of partitions of $39$
  which dominate $(29,6,4)$ and which do not have empty $3$-core. Then $\Lambda$
  is $3$-cosaturated and contains the 10 partitions listed below. The reader may
  check that each of these partitions satisfies the equivalent conditions of
  \autoref{projective}. Consequently, the decomposition matrix of
  $\S_{\C,\omega}(\Lambda)$, where $\omega=\exp(2\pi i/3)$, is the identity
  matrix and $\S_{\C,\omega}(\Lambda)$ is semisimple. In contrast, if
  $(e,p)=(3,2)$ then, using \cite{SPECHT} together with the Steinberg tensor
  product theorem via~\cite[Proposition~5.4]{Cox:blocks} (see section~3.5),
 the decomposition matrix of~$\S_{k,q}(\Lambda)$ is the following.
   $$
   \begin{array}{l|*{10}c}
      (34,5)   &1& & & & & & & & &\\
      (31,5,3) &.&1& & & & & & & &\\
      (31,8)   &.&1&1& & & & & & &\\
      (37,2)   &.&1&.&1& & & & & &\\
      (30,7,2) &.&.&.&.&1& \\
      (33,4,2) &.&.&.&.&.&1\\
      (29,6,4) &.&.&.&.&.&.&1\\
      (29,9,1) &.&.&.&.&.&.&1&1\\
      (32,6,1) &.&.&.&.&.&.&1&1&1\\
      (35,3,1) &.&.&.&.&.&.&.&.&1&1\\
   \end{array}
   $$
   Therefore, when $p=2$ the partitions in $\Lambda$ which satisfy the
   equivalent conditions of \autoref{projective} are $(34,5)$, $(31,8)$,
   $(37,2)$, $(30,7,2)$, $(33,4,2)$ and $(35,3,1)$.
\end{Example}

The hook $R^\lambda_{ab}$ is \textbf{horizontal} if it is entirely
contained in row~$a$ of $\lambda$. Thus, $R^\lambda_{ab}$ is horizontal if
and only if~$\lambda_b'=a$.

\begin{Definition}\label{D:horizontal}
  Suppose that $\lambda\in\Lambda$. Then $\lambda$ \textbf{only contains horizontal
  hooks} if $\lambda$ is either an $e$-core or every removable $e$-hook of~$\lambda$
  is horizontal and unwrapping any of these hooks gives
  a partition which only contains horizontal $e$-hooks.
\end{Definition}

\begin{Example}
  Suppose that $e=3$ and $\mu=(7,4)$. Then $\mu$ does \textit{not} contain only
  horizontal $e$-hooks because, even though all of the removable $3$-hooks of
  $\mu$ are horizontal, $R^{(7,4)}_{15}$, $R^{(4^2)}_{13}$, $R^{(3,2)}_{12}$ is
  a sequence of $e$-hooks leading to its $3$-core $(1^2)$, however, only the
  first of these hooks is horizontal.
\end{Example}

We now come to the main result of this section, which is both tricky to prove
and pivotal for our Main Theorem.  In particular, \autoref{conditions}
shows that all of the partitions in \autoref{first example} contain only
horizontal $3$-hooks.

\begin{Proposition}\label{conditions}
  Suppose that $\mu\in\Lambda$ and let $\kappa$
  be the $e$-core of~$\mu$. Then the following are equivalent:
  \begin{enumerate}
    \item If $(a,c),(b,c)\in\Diag(\mu)$ then $e$ divides $h^\mu_{ac}$ if and
      only if $e$ divides $h^\mu_{bc}$.
    \item The partition $\mu$ only contains horizontal $e$-hooks.
    \item $\mu_i-\mu_{i+1}\equiv-1\pmod e$, whenever $i\ge1$ and
      $\mu_{i+1}>\kappa_{i+1}$.
  \end{enumerate}
  Moreover, these three combinatorial conditions are all equivalent to
  $\Delta^\mu_k=P^\mu_k$ being projective as an $\S_{k,q}(\Lambda)$-module when
  $p=0$.
\end{Proposition}

\begin{proof}
  By \autoref{projective}, part~(a) is equivalent to $\Delta^\mu_k$
  being projective when $p=0$,  so it is enough to show that (a)--(c) are
  equivalent. Before we start, recall that a partition $\nu$ is an $e$-core if
  and only if $e$ does not divide $h^\nu_{ab}$, for all $(a,b)\in\Diag(\nu)$.
  This is easily proved using \autoref{wrapping}.

  Let $w$ be the $e$-weight of~$\mu$. If $w=0$ then $\mu$ is an $e$-core and
  (a)--(c) are equivalent because all three statements are vacuous by the
  remarks in the first paragraph. We now assume that $w>0$ and argue by
  induction on~$w$. Suppose that~(a) holds and fix $(a,c)\in\Diag(\mu)$ where
  $a$ is maximal such that $e$ divides $h^\mu_{ac}$. By part~(a), and the maximality
  of $a$, the node $(a,c)$ must be at the bottom of its column.  Hence,
  $R^\mu_{ac}$ is a horizontal hook and, by changing~$c$ if necessary, we may
  assume that $h^\mu_{ac}=e$. Let $\nu$ be the partition obtained from~$\mu$ by
  unwrapping $R^\mu_{ac}$. Since~$R^\mu_{ac}$ is horizontal, and $h^\mu_{ac}=e$,
   $$h^\nu_{xy}=\begin{cases}
    h^\mu_{xy}-1,&\text{if $c\le y<c+e$ and $x<a$},\\
    h^\mu_{xy}-e,&\text{if $y<c$ and $x\le a$},\\
    h^\mu_{xy},&\text{otherwise.}
  \end{cases}$$
  Therefore, $\nu$ satisfies~the condition in part~(a) of the Proposition.
  Hence, by induction on~$w$,~$\nu$ also satisfies condition~(b) so that every
  removable $e$-hook contained in~$\nu$ is horizontal.  Now suppose that
  $R^\mu_{xy}$ is any removable $fe$-hook in ~$\mu$, for $f\ge1$. If
  $R^\mu_{xy}\cap R^\mu_{ac}=\emptyset$ then
  $R^\mu_{xy}\subseteq[\nu]$ so that it is a union of horizontal $e$-hooks by
  induction. If $R^\mu_{xy}\cap R^\mu_{ac}$ is non-empty then, $y\notin(c,c+e)$
  by~(a) since $h^\mu_{xy}=fe$. Therefore, $R^\mu_{ac}\subseteq R^\mu_{xy}$, so
  that $R^\mu_{xy}{\setminus}R^\mu_{ac}$ is a union of horizontal $e$-hooks
  in~$\nu$. Continuing in this way shows that every $e$-hook contained in~$\mu$
  is either equal to $R^\mu_{ac}$ or it is an $e$-hook contained in~$\nu$.
  Therefore, by induction, all of the $e$-hooks contained in~$\mu$ are
  horizontal so that~(b) holds.

  Now suppose that (b) holds. By way of contradiction, suppose that
  $\mu_i-\mu_{i+1}\not\equiv-1\pmod e$, for some $i$ with
  $\mu_{i+1}>\kappa_{i+1}$. Without loss of generality, we may assume that $i$
  is maximal with this property. Let~$c$ be the unique integer such that $0<c<e$
  and $e-c-1\equiv\mu_i-\mu_{i+1}$. Since $i$ was chosen to be maximal,
  $(i+1,c')$ is at the bottom of its column, where $c'=\mu_{i+1}-c+1$, so that
  $R^\mu_{ic'}$ is a removable $fe$-hook which is not horizontal, for some
  $f\ge1$. It follows that $\mu$ contains a non-horizontal $e$-hook, which is a
  contradiction. Hence, $\mu_i-\mu_{i+1}\equiv-1\pmod e$ whenever
  $\mu_{i+1}>\kappa_{i+1}$ and (c) holds.

  Finally, suppose that (c) holds. Suppose that $i$ is maximal such that
  $\mu_{i+1}>\kappa_{i+1}$. Now because $\mu_i-\mu_{i+1}\equiv-1\pmod e$,
  row $i+1$ of $\mu$ contains $(\mu_{i+1}-\kappa_{i+1})/e$ horizontal
  $e$-hooks and, moreover, $\mu_i>\kappa_i$. Hence, removing these
  horizontal $e$-hooks and arguing by induction it follows that
  $\mu_j-\mu_{j+1}\equiv-1\pmod e$, for $1\le j\le i$.

  Now let $(a,c)$ and $(b,c)$ be two nodes in $\Diag(\mu)$ with $a<b$. If
  $(b,c)\in[\kappa]$ then $(a,c)\in[\kappa]$ so that $e$ does not divide
  $h^\mu_{ac}$ and $e$ does not divide $h^\mu_{bc}$. If
  $(b,c)\notin[\kappa]$ then $\mu_b>\kappa_b$ and, by the last paragraph,
  $\mu_a>\kappa_a$ since $a<b$. Let $\mu'=(\mu'_1,\mu'_2,\dots)$ be the partition
  which is conjugate to $\mu$. Then,
  \begin{align*}
    h^\mu_{ac}-h^\mu_{bc}
       &=(\mu_a-a+\mu_c'-c+1)-(\mu_b-b+\mu'_c-c+1)\\
       &=\mu_a-\mu_b+b-a
        \equiv 0\pmod e,
  \end{align*}
  where the last congruence follows because
  $\mu_j-\mu_{j+1}\equiv-1\pmod e$, for $1\le j\le i$. Hence, (a) holds and the proof is
  complete.
\end{proof}

\subsection{The combinatorics of our Main Theorem}
Using \autoref{conditions} we can now properly define the combinatorics
underpinning our Main Theorem. The main result of this section is
\autoref{Frobenius} which gives a combinatorial reduction of the calculation of the
Jantzen coefficients to the case when $s_\Lambda(\mu)=1$.

Recall from the introduction that $\P=\{1,e,ep,ep^2,\dots\}$. Suppose that
$\mu\in\Lambda$ has $e$-core $\kappa=(\kappa_1,\kappa_2,\dots)$.  Let
$\Lambda_\kappa$ be the set of partitions in $\Lambda$ which have
$e$-core~$\kappa$ and define the \textbf{length function}
$\ell_\Lambda\map\Lambda\N$ by
$$\ell_\Lambda(\mu)=\min\set{i\ge0|\lambda_j=\kappa_j \text{ whenever }j>i
\text{ and } \lambda\in\Lambda_\kappa}.$$
(By definition a partition is an infinite non-increasing sequence
$\mu=(\mu_1,\mu_2,\dots)$ so this makes sense.) Observe that if
$\Lambda\subseteq\Lambda_r$ then $\ell_\Lambda(\mu)<r$, for all $\mu\in\Lambda$.
Moreover, $\ell_\Lambda(\mu)=0$ if and only if $\mu=\kappa$ is an $e$-core and
$\ell_\Lambda(\mu)=1$ only if $\Lambda_\kappa=\{\mu\}$.

The reason why the length function $\ell_\Lambda$ is important is that
if $\kappa$ is an $e$-core and if~$\mu$ is any partition in
$\Lambda_\kappa$ then $\mu_i=\kappa_i$, whenever $i>\ell_\Lambda(\mu)$. In
particular, when applying the sum formula we can never move $e$-hooks below
row~$\ell_\Lambda(\mu)$.

Following the introduction, if $\mu\in\Lambda$ and $\ell_\Lambda(\mu)\le1$ set
$s_\Lambda(\mu)=1$ and otherwise define
$$s_\Lambda(\mu)=\max\set{s\in\P|
                    \mu_i-\mu_{i+1}\equiv-1 \pmod s, 1\le i<\ell_\Lambda(\mu)}.
$$
This definition is stronger than it appears.

\begin{Lemma}\label{chi}
    Suppose that $\mu\in\Lambda$ and that $s'\in\P$ with $0<s'\le s_\Lambda(\mu)$. Then
    $$\mu_i-\mu_{i+1}\equiv -1\pmod{s'},$$
    for $1\le i<\ell_\Lambda(\mu)$. Moreover, every removable $s'$-hook contained
    in~$\mu$ is horizontal and if~$(a,c),(b,c)\in\Diag(\mu)$ then $s'$ divides
    $h^\mu_{ac}$ if and only if $s'$ divides $h^\mu_{bc}$.
\end{Lemma}

\begin{proof}
  If $0\ne s'\in\P$ and $s_\Lambda(\mu)\ge s'$ then $\mu_i-\mu_{i+1}\equiv-1\pmod{s'}$
  since $s'$ divides $s$. Applying \autoref{conditions} with $e=s'$ shows that
  every removable $s'$-hook contained in $\mu$ is horizontal and that $s'$
  divides $h^\mu_{ac}$ if and only if $s'$ divides $h^\mu_{bc}$, for
  $(a,c),(b,c)\in\Diag(\mu)$.
\end{proof}

Armed with \autoref{chi} we can now give a more transparent definition of the partition
$\chi_\Lambda(\mu)=(\chi_1,\chi_2,\dots)$ from the introduction. That is, let
$\chi_i$ be the number of horizontal $s$-hooks in row~$i$ of $\mu$,
where $s=s_\Lambda(\mu)$ and $1\le i\le\ell_\Lambda(\mu)$. Hence, if $\kappa^{(s)}$ is
the $s$-core of~$\mu$ then
$$\mu_i\equiv\kappa^{(s)}_i\pmod{s}, \quad\text{for } 1\le i\le\ell_\Lambda(\mu).$$
Therefore, since $\mu$ is a partition which contains only horizontal $s$-hooks, the
number of $s$-hooks in row~$i$ of~$\mu$ is greater than that number of $s$-hooks in
row~$i+1$. Hence, $\chi_\Lambda(\mu)$ is a partition.

We illustrate all of the definitions in this section in \autoref{Ex:DNF}
below.


By definition, if $\lambda$ and $\mu$ are two partitions in~$\Lambda$ which have
the same $e$-core then $\ell_\Lambda(\lambda)=\ell_\Lambda(\mu)$. This
observation accounts for the dependence of the integer $s_\Lambda(\mu)$ and the
partition $\chi_\Lambda(\mu)$ upon the poset~$\Lambda$.

\begin{Definition}\label{DNF}
  Define $\sim_\Lambda$ to be the equivalence relation on $\Lambda$ such
  that $\lambda\sim_\Lambda\mu$, for $\lambda,\mu\in\Lambda$, if
  \begin{enumerate}
    \item $\lambda$ and $\mu$ have the same $e$-core;
    \item $s_\Lambda(\lambda)=s_\Lambda(\mu)$; and,
    \item if $s_\Lambda(\mu)>1$ then $\chi_\Lambda(\lambda)$ and $\chi_\lambda(\mu)$
    have the same $p$-core.
  \end{enumerate}
\end{Definition}

For part~(c), recall that if $p=0$ then the $0$-core of the partition $\nu$ is $\nu$.

Thus, our Main Theorem says that if $\lambda,\mu\in\Lambda$ then
$\Delta^\lambda_k$ and $\Delta^\mu_k$ are in the same block if and only if
$\lambda\sim_\Lambda\mu$. We prove this in the next section. First,
however, we give an example and begin to investigate the combinatorics of
the equivalence relation $\sim_\Lambda$.

\begin{Example}\label{Ex:DNF}
    Suppose that $(e,p)=(3,2)$ and let
    $\Lambda$ be the set of partitions of $39$ which dominate $(29,6,4)$
    and which do not have empty $3$-core. Then $\Lambda$ is $3$-cosaturated
    and~it contains the $10$ partitions in the table below which
    describes the equivalence $\sim_\Lambda$.
$$\begin{array}{l|lc|cll}
  \multicolumn1c\mu & \core[3](\mu) & \ell_\Lambda(\mu) & s_\Lambda(\mu)
         & \chi_\Lambda(\mu) & \core[2](\chi_\Lambda(\mu))\\\toprule
(35,3,1)&(5,3,1)&3& 3&(10) &(0)\\
(32,6,1)&(5,3,1)&3& 3&(9,1) &(0)\\
(29,9,1)&(5,3,1)&3& 3&(8,2)  &(0)\\
(29,6,4)&(5,3,1)&3& 3&(8,1,1)&(0)\\\midrule
(33,4,2)&(6,4,2)&2& 6&(4)&(0)\\\midrule
(30,7,2)&(6,4,2)&2&24&(0)&(0)\\\midrule
(37,2  )&(4,2)  &3& 3&(11)&(1)\\
(31,8  )&(4,2)  &3& 3&(9,2)&(1)\\
(31,5,3)&(4,2)  &3& 3&(9,1,1)&(1)\\\midrule
(34,5  )&(4,2)  &3& 6&(4)&(0)\\\bottomrule
\end{array}$$
The different regions in the table give the $\sim_\Lambda$ equivalence classes
in~$\Lambda$. By our Main Theorem these regions label the blocks of
$\S_{k,q}(\Lambda)$. The reader can check that this agrees with
\autoref{first example} which gives the decomposition matrix for
$\S_{k,q}(\Lambda)$. This example shows that $\ell_\Lambda$ need not be constant
on $\Lambda$.

Continuing \autoref{first example}, if $(e,p)=(3,0)$ then $\ell_\Lambda(\mu)$ is
as given above but $s_\Lambda(\mu)=3$ for all $\mu\in\Lambda$. Therefore, by our Main
Theorem, all of these partitions are in different blocks because the partitions
$\core[0]\!\big(\chi_\Lambda(\mu)\big)=\chi_\Lambda(\mu)$ are distinct, for $\mu\in\Lambda$.
Once again, this agrees with the block decomposition of $\S_{k,q}(\Lambda)$
given in \autoref{first example} when $(e,p)=(3,0)$.
\end{Example}

The following results establish properties of the equivalence relation
$\sim_\Lambda$ that we need to prove our Main Theorem.

\begin{Lemma}\label{L:sRunners}
  Suppose that $\mu\in\Lambda$ and $s\in\P$. Then $s_\Lambda(\mu)\ge s$ if and
  only if the last $\ell_\Lambda(\mu)$ beads on an $s$-abacus configuration for
  $\mu$ are all on the same runner.
\end{Lemma}

\begin{proof}
  Observe that  $\mu_i-\mu_{i+1}\equiv-1\pmod{s}$, for $1\le i<\ell$, if and only if
  $$ \mu_{i+1}-(i+1)\equiv\mu_i-i\pmod{s}, $$
  whenever $1\le i<\ell$. For any positive integer $m$ the $m$-beta numbers for $\mu$ are
  $m+\mu_j-j$, for $1\le j\le m$. Hence, using \autoref{chi}, it follows that
  $s_\Lambda(\mu)\ge s$ if and only if the last~$\ell$ beads on any abacus configuration
  for~$\mu$ with~$s$ runners all lie on the same runner.
\end{proof}

\begin{Lemma}\label{chi hooks}
  Suppose that $\mu\in\Lambda$ and that $(a,b)\in[\chi]$ is a node in
  $\chi=\chi_\Lambda(\mu)$ and let $s=s_\Lambda(\mu)$. Then
  $h^\chi_{ab}=\tfrac1sh^\mu_{aB},$
  where column~$B$ of $\mu$ is the $b^{\text{th}}$ column of
  $\mu$ with hook lengths divisible by~$s$, reading from left to right.
\end{Lemma}

\begin{proof}
  By \autoref{conditions} and the definition of $s=s_\Lambda(\mu)$, all
  of the removable $s$-hooks in $\mu$ are horizontal so the definition
  of~$B$ makes sense. The Lemma follows from the observation that the nodes
  in~$\chi$ correspond to the removable $s$-hooks in~$\mu$ and that the
  nodes in the rim rook $R^\chi_{ab}$ correspond to the removable
  $s$--hooks which make up the rim hook $R^\mu_{aB}$ in~$\mu$.
\end{proof}

\begin{Lemma}\label{chi dominance}
  Suppose that $\lambda$ and $\mu$ have the same $e$-core and
  that $s_\Lambda(\lambda)=s_\Lambda(\mu)$, for partitions $\lambda,\mu\in\Lambda$. Then
  $\lambda\gedom\mu$ if and only if
  $\chi_\Lambda(\lambda)\gedom\chi_\Lambda(\mu)$.
\end{Lemma}

\begin{proof}By \autoref{conditions} the partitions $\lambda$ and $\mu$ only
  contain horizontal $s$-hooks, where $s=s_\Lambda(\mu)$. The Lemma follows using
  this observation and the correspondence between the nodes in
  $\chi_\Lambda(\lambda)$ and $\chi_\Lambda(\mu)$ and the horizontal
  $s$-hooks in $\lambda$ and $\mu$, respectively.
\end{proof}

The following result is a key reduction step for understanding the blocks of
$\S_{k,q}(\Lambda)$. \autoref{Frobenius} can be interpreted as saying that the
Steinberg tensor product theorem preserves Jantzen equivalence --- note,
however, that its proof requires no knowledge of Steinberg.  This result will
allow us to reduce Jantzen equivalence to the case where $s_\Lambda(\mu)=1$.
Recall that $\S_{k,1}(\Lambda)$ is the Schur algebra with parameter $q=1$.

\begin{Proposition}\label{Frobenius}
  Suppose that $\lambda,\mu\in\Lambda_\kappa$, where~$\kappa$ is an $e$-core,
  are partitions with $s=s_\Lambda(\lambda)=s_\Lambda(\mu)>1$. Let
  $\Gamma=\set{\chi_\Lambda(\nu)|\nu\in\Lambda_\kappa\text{ and } s_\Lambda(\nu)=s}.$
  Then $\Gamma$ is an cosaturated set of partitions and
  $$J^\Lambda_{\lambda\mu}\ne0 \quad\text{if and only if}\quad
  J^\Gamma_{\chi_\Lambda(\lambda)\chi_\Lambda(\mu)}\ne0,$$ where
  $J^\Gamma_{\chi_\Lambda(\lambda)\chi_\Lambda(\mu)}$ is a Jantzen coefficient
  for the algebra $\S_{k,1}(\Gamma)$. Moreover, if $p=0$ then
  $J_{\lambda\mu}^\Lambda=0$.
\end{Proposition}

\begin{proof}
  By \autoref{chi dominance}, $\Gamma$ is an $e$-cosaturated set of partitions. In
  order to compare the Jantzen coefficients of the algebras $\S_{k,q}(\Lambda)$ and
  $\S_{k,1}(\Gamma)$ write $s=ep^d$, for some $d\ge0$ (with $d=0$ if $p=0$). By
  \autoref{chi}, if $(x,z)\in\Diag(\mu)$ then $\nu_{e,p}(h^\mu_{xz})\ne0$ only if $s$
  divides $h^\mu_{xz}$. Moreover, using the notation of \autoref{chi hooks}, if $s$
  divides $h^\mu_{aB}$, for $(a,B)\in\Diag(\mu)$, then
  $$\nu_{e,p}(h^\mu_{aB}) = \nu_{e,p}(sh^{\chi_\Lambda(\mu)}_{ab})
            =\nu_{e,p}(ep^dh^{\chi_\Lambda(\mu)}_{ab})
            =\begin{cases}
              d+\nu_{p,p}(h^{\chi_\Lambda(\mu)}_{ab}),&\text{if }p>0,\\
              1&\text{if }p=0.
            \end{cases}
  $$
  Hence,  $J^\Lambda_{\lambda\mu}\ne0$ if and only
  if~$J^\Gamma_{\chi_\Lambda(\lambda)\chi_\Lambda(\mu)}\ne0$ by
  \autoref{non-vanishing cols}. Finally, if $p=0$ then
  $J_{\lambda\mu}^\Lambda=0$ by \autoref{non-vanishing cols} because, by
  what we have shown, if $s_\Lambda(\mu)>1$ then $\nu_{e,0}(h^\mu_{xy})$ is constant
  on the columns of~$\mu$.
\end{proof}

\subsection{The Main Theorem}
We are now ready to prove our main theorem. We start by settling the case when
$s_\Lambda(\mu)=1$.

Recall from after \autoref{D:Cosaturated} that if $\kappa$ is an $e$-core then
$\Lambda_\kappa$ is the set of partitions in $\Lambda$ with $e$-core~$\kappa$.
By~\eqref{E:SLambdakappa} the algebra $\S_{k,q}(\Lambda_\kappa)$ is a direct summand
of~$\S_{k,q}(\Lambda)$.

\begin{Lemma}\label{s=1}
  Suppose that $s_\Lambda(\tau)=1$, where $\tau\in\Lambda$ has $e$-core $\kappa$. Then
  $$\set{\mu\in\Lambda|\mu\sim_\Lambda\tau}=\Lambda_\kappa
              =\set{\mu\in\Lambda|\mu\sim_J\tau}.$$
  In particular, $\mu\sim_\Lambda\tau$ if and only if $\mu\sim_J\tau$.
\end{Lemma}

\begin{proof}
  By definition, if $\mu\in\Lambda$ then $\mu\sim_\Lambda\tau$ if and only if
  $\mu\in\Lambda_\kappa$ and $s_\Lambda(\mu)=1$. Suppose, by way of
  contradiction, that $s_\Lambda(\mu)>1$ for some $\mu\in\Lambda_\kappa$.
  Taking $s'=e$ in \autoref{chi}, it follows that
  $$-1\equiv\mu_i-\mu_{i+1}\pmod e,\qquad\text{for }1\le i<\ell_\Lambda(\mu).$$
  Combining parts~(b) and (c) of \autoref{conditions}, this last
  equation is equivalent to
  $$-1\equiv\kappa_i-\kappa_{i+1}\pmod e,\qquad\text{for }1\le i<\ell_\Lambda(\mu).$$
  By the same argument, since $\ell_\Lambda(\tau)=\ell_\Lambda(\mu)$, this implies that
  $$-1\equiv\tau_i-\tau_{i+1}\pmod e,\qquad\text{for }1\le i<\ell_\Lambda(\tau).$$
  This implies that $s_\Lambda(\tau)\ge e$, a contradiction! Therefore,
  $s_\Lambda(\mu)=1$, for all $\mu\in\Lambda_\kappa$. Hence,
  $\Lambda_\kappa=\set{\mu\in\Lambda|\mu\sim_\Lambda\tau}$, giving the left
  hand equality of the Lemma.

  We now show that $\mu\sim_J\tau$ if and only if $\mu\in\Lambda_\kappa$.  If
  $\mu\sim_J\tau$ then $\mu\in\Lambda_\kappa$ by \autoref{e-cores}. To prove
  the converse, let $\gamma$ be the unique partition with $e$-core $\kappa$
  which has $(|\tau|-|\kappa|)/e$ horizontal $e$-hooks in its first row. Then
  $\gamma\gedom\mu$ for all $\mu\in\Lambda_\kappa$. To complete the proof it is
  enough to show that $\mu\sim_J\gamma$, whenever $\mu\in\Lambda_\kappa$. If
  $\mu=\gamma$ there is nothing to prove, so suppose that $\mu\ne\gamma$. If
  $J_{\lambda\mu}\ne0$ for some $\lambda\in\Lambda$ then $\lambda\gdom\mu$ by
  \autoref{e-cores} so that $\mu\sim_J\lambda\sim_J\gamma$ by induction on
  dominance.

  We have now reduced to the case when $J_{\lambda\mu}=0$ for all
  $\lambda\in\Lambda$. Consequently,
  $\nu_{e,p}(h^\mu_{ac})=\nu_{e,p}(h^\mu_{bc})$, for all
  $(a,c),(b,c)\in\Diag(\mu)$ by \autoref{projective}(e). Hence, $\mu$
  contains only horizontal $e$-hooks by \autoref{conditions}. On the
  other hand, since $s_\Lambda(\mu)=1$ by the last paragraph, there exists an
  integer $i$ such that $\mu_i-\mu_{i+1}\not\equiv-1\pmod e$ and $1\le
  i<\ell_\Lambda(\mu)$. Fix $i$ which is minimal with this property and notice
  that we must have $\mu_{i+1}=\kappa_{i+1}$ by \autoref{conditions}(c).
  Recalling that all of the $e$-hooks in $\mu$ are horizontal, let $\lambda$ be
  the partition obtained by unwrapping the lowest removable $e$-hook from $\mu$
  and then wrapping it back on with its foot node in row $i+1$. Then $\lambda$
  is a partition because~$i$ is minimal such that
  $\mu_i-\mu_{i+1}\not\equiv-1\pmod e$. Moreover, since $\Lambda$ is
  cosaturated, $\lambda\in\Lambda_\kappa$ because $\ell_\Lambda(\mu)>i$,
  $\mu_{i+1}=\kappa_{i+1}$ and all of the $e$-hooks in~$\mu$ are horizontal.
  Next observe that $J_{\mu\lambda}\ne0$ by \autoref{non-vanishing rows}
  because the valuations of the corresponding hook lengths are different since
  all of the $e$-hooks in~$\mu$ are horizontal. Now let $\sigma$ be the
  partition obtained by unwrapping this same hook from $\lambda$ and wrapping it
  back on as a horizontal hook in the first row. By construction all of the $e$-hooks in
  $\sigma$ are horizontal so, as before, $J_{\sigma\lambda}\ne0$
  by \autoref{non-vanishing cols} (note, however, that
  $J_{\sigma\mu}=0$). Hence, $\mu\sim_J\lambda\sim_J\sigma\sim_J\gamma$, with
  the last equivalence following by induction since $\sigma\gdom\mu$. This completes
  the proof.
\end{proof}

\begin{Lemma}\label{J=>S}
  Suppose that $\lambda\sim_J\mu$, for $\lambda,\mu\in\Lambda$. Then
  $s_\Lambda(\lambda)=s_\Lambda(\mu)$.
\end{Lemma}

\begin{proof} Let $s=s_\Lambda(\mu)$ and let $\ell=\ell_\Lambda(\mu)$. By
  \autoref{s=1} we may assume that $s>1$ and hence that $p>0$ since
  $s_\Lambda(\mu)\in\{1,e\}$ if $p=0$.  It is enough to show that
  $s_\Lambda(\lambda)=s$ whenever $J_{\lambda\mu}\ne0$.  By
  \autoref{non-vanishing cols}, $J_{\lambda\mu}\ne0$ if and only if there
  exist nodes $(x,z),(y,z)\in\Diag(\mu)$ such that $x<y\le\ell$,
  $\nu_{e,p}(h^\mu_{xz})\ne\nu_{e,p}(h^\mu_{yz})$ and $\lambda$ is obtained from
  $\mu$ by unwrapping~$R^\mu_{yz}$ and wrapping it back on with its foot node in
  row~$x$.  Therefore,~$s$ divides both of~$h^\mu_{xz}$ and~$h^\mu_{yz}$ and
  $\lambda$ is obtained from $\mu$ by moving a union of $s$-hooks.

  By \autoref{L:sRunners} the last~$\ell$ beads are always on the same runner
  in any $s$-abacus. Therefore, by the last paragraph, an $s$-abacus for $\lambda$ is
  obtained from the $s$-abacus configuration for~$\mu$ by moving two beads on
  the same runner. That is, the abacus configuration for $\lambda$ is obtained
  from an $s$-abacus for $\mu$ by moving one bead up $fs$ positions and another
  bead down $fs$-positions, for some $f\ge1$. Hence,
  $s_\Lambda(\lambda)\ge s_\Lambda(\mu)=s$
  by \autoref{L:sRunners} (since $\ell_\Lambda(\lambda)=\ell$ by
  \autoref{e-cores}).

  By symmetry, using \autoref{non-vanishing rows}
  instead of \autoref{non-vanishing cols},
  $s_\Lambda(\mu)\ge s_\Lambda(\lambda)$. Hence, $s_\Lambda(\mu)=s_\Lambda(\lambda)$
  as required.
\end{proof}

We can now prove our Main Theorem.

\begin{proof}[Proof of the Main Theorem]
  By \autoref{blocks} we need to prove that $\lambda\sim_J\mu$ if
  and only if $\lambda\sim_\Lambda\mu$, for $\lambda,\mu\in\Lambda$. If $p=0$ then the result
  follows from \autoref{Frobenius} and \autoref{s=1}, so assume that $p>0$.

  First suppose that $\lambda\sim_\Lambda\mu$, for $\lambda,\mu\in\Lambda$.
  To show that $\lambda\sim_J\mu$ we argue by induction on $s=s_\Lambda(\mu)$.
  If $s=1$ the result is just \autoref{s=1}, so suppose that $s>1$. As in
  \autoref{Frobenius} let
  $\Gamma=\set{\chi_\Lambda(\nu)|\nu\in\Lambda_\kappa\text{ and }
                                    s_\Lambda(\nu)=s}$,
  an $e$-cosaturated set of partitions. By definition,
  $s_\Gamma(\chi_\Lambda(\lambda))=1=s_\Gamma(\chi_\Lambda(\mu))$ and, since
  $\lambda\sim_\Lambda\mu$, the partitions $\chi_\Lambda(\lambda)$
  and~$\chi_\Lambda(\mu)$ have the same $p$-core. Therefore,
  $\chi_\Lambda(\lambda)\sim_{J^\Gamma}\chi_\Lambda(\mu)$ by \autoref{s=1}.
  Hence, by \autoref{Frobenius}, $\lambda\sim_J\mu$  as required.

  To prove the converse it is enough to show that $\lambda\sim_\Lambda\mu$
  whenever $J_{\lambda\mu}\ne0$. By \autoref{e-cores}, $\lambda$ and $\mu$
  have the same $e$-core. Moreover, $s_\Lambda(\lambda)=s_\Lambda(\mu)$, by
  \autoref{J=>S}. Finally,
  $\chi_\Lambda(\lambda)\sim_{J^\Gamma}\chi_\Lambda(\mu)$ are Jantzen equivalent
  for $\S_{k,1}(\Gamma)$ by \autoref{Frobenius} since
  $\lambda\sim_J\mu$. Consequently, $\chi_\Lambda(\lambda)$ and
  $\chi_\Lambda(\mu)$ have the same $p$-core by \autoref{e-cores}. Hence,
  $\lambda\sim_\Lambda\mu$ as we wanted to show.
\end{proof}

\section*{Acknowledgement}
We thank Steve Donkin and Hebing Rui for their helpful comments on preliminary
versions of this paper.



\end{document}